\documentclass[12pt]{amsart}

\title{Separated Lie models and the homotopy Lie algebra}
\author{Peter Bubenik}
\address{Department of Mathematics, Cleveland State University, 2121 Euclid Ave. RT 1515, Cleveland OH, 44115--2214, USA}
\email{p.bubenik@csuohio.edu}
\urladdr{http://academic.csuohio.edu/bubenik\_p/}
\date{May 7, 2007}

\usepackage{amsmath}
\usepackage{amsthm}
\usepackage{amssymb}
\usepackage[all]{xy}
\usepackage{hyperref}

\newtheorem{theorem}{Theorem}[section]
\newtheorem{lemma}[theorem]{Lemma}
\newtheorem{prop}[theorem]{Proposition}

\newtheorem{corollary}[theorem]{Corollary}

\theoremstyle{definition}
\newtheorem{defn}[theorem]{Definition}

\theoremstyle{remark}
\newtheorem{rem}[theorem]{Remark}
\newtheorem{example}[theorem]{Example}

\numberwithin{equation}{section}

\newcommand {\F} {\ensuremath {\mathbb{F}} }

\newcommand {\Q} {\ensuremath {\mathbb{Q}} }

\newcommand {\kk} {\ensuremath {\mathbf{k}} }
\newcommand {\CP} {\ensuremath {\mathbb{CP}} }

\newcommand {\lX} {\ensuremath {\Omega X} }

\newcommand {\eL} {\mathbf{\underline{L}}} 
\newcommand {\freeL} {\ensuremath {\mathbb{L}} }
\newcommand {\freeLab} {\ensuremath {\mathbb{L}_{ab}} }
\newcommand {\isom} {\ensuremath {\cong} }
\newcommand {\tensor} {\ensuremath {\otimes} }
\newcommand {\incl} {\ensuremath {\hookrightarrow} }
\newcommand {\injects} {\ensuremath {\hookrightarrow} }

\newcommand {\rsdp} {\ensuremath {\rtimes} }

\newcommand {\onto} {\ensuremath {\twoheadrightarrow} }
\newcommand {\isomto} {\ensuremath {\xrightarrow{\isom}} }

\newcommand {\xto}[1] {\ensuremath {\xrightarrow{#1}} }

\newcommand {\clQ} {\cl_{\Q}}

\newcommand {\rht} {\simeq_{\Q}}
\newcommand {\catQ} {\cat_{\Q}}
\newcommand {\qi} {\xto{\simeq}}

\DeclareMathOperator{\im}{im}
\DeclareMathOperator{\gr}{gr}
\DeclareMathOperator{\id}{id}

\DeclareMathOperator{\cl}{cl}

\DeclareMathOperator{\Hom}{Hom}
\DeclareMathOperator{\cat}{cat}

\begin{document}

\begin{abstract}
  A simply connected topological space $X$ has homotopy Lie algebra $\pi_*(\lX) \tensor \Q$. Following Quillen, there is a connected differential graded free Lie algebra (dgL) called a Lie model, which determines the rational homotopy type of $X$, and whose homology is isomorphic to the homotopy Lie algebra. We show that such a Lie model can be replaced with one that has a special property we call \emph{separated}. The homology of a separated dgL has a particular form which lends itself to calculations.
%
\end{abstract}

\maketitle

\section{Introduction}

All of our topological spaces are assumed to be simply connected, and have finite rational homology in each dimension. For such a space $X$, there exists a differential graded Lie algebra $(\freeL V,d)$, where $\freeL V$ denotes the free graded Lie algebra on the rational vector space $V$, called a \emph{Lie model}, which determines the rational homotopy type of $X$, and whose homology is isomorphic to $\pi_*(\lX) \tensor \Q$, the \emph{(rational) homotopy Lie algebra} of $X$. While this description is pleasing in principle, it is less satisfying when explicit calculations are desired. Here, we want to put a certain structure on a Lie model that will prove conducive to calculation.

Assume that our Lie algebras are finite in each dimension and concentrated in positive dimension. Our assumption on our spaces implies that their Lie models satisfy this assumption. We are particularly interested in dgLs where $V$ can be bigraded as follows: $V_* = \oplus_{i=1}^N V_{i,*}$, such that 
\begin{equation} \label{eq:condL}
dV_{i,*} \subset \freeL (\oplus_{j=1}^{i-1} V_{j,*}).
\end{equation}
We will call the first gradation \emph{degree} and the second, usual gradation \emph{dimension}. 

Our \emph{separated} condition will say that the attaching map from degree $i+1$ has no effect on the homology coming from degree $\leq i-1$. We need to introduce some notation to make this precise.

Given a Lie algebra $L = (\freeL V,d)$ as above, let 
\begin{equation} \label{eq:Li}
L_i = (\freeL V_{\leq i,*}, d). 
\end{equation}
Let $ZL_i$ and $BL_i$ denote the cycles and boundaries in $L_i$. There is an induced map 
\begin{equation} \label{eq:tilded}
\tilde{d}_i : V_{i,*} \xto{d} Z(L_{i-1}) \onto HL_{i-1}.
\end{equation}
The inclusion $L_i \incl L_{i+1}$ induces inclusions $BL_i \incl
BL_{i+1}$ and $ZL_i \incl ZL_{i+1}$.
Thus there is an induced map $HL_i \to HL_{i+1}$.

We introduce some notation for two Lie subalgebras of $HL_i$ to which we will often refer. Let $HL_i^-$ be the Lie subalgebra of $HL_i$
given by $\im (HL_{i-1} \to HL_i)$, and let $HL_i^+$ be the Lie ideal generated by $\im(\tilde{d}_{i+1}: V_{i+1,*} \to HL_i)$.

\begin{defn} \label{def:separated}
Let $L$ be a dgL satisfying~\eqref{eq:condL}. Say that $L$ is \emph{separated} if for all $i$,
\[ HL_i^+ \cap HL_i^- = 0. \]
\end{defn}

Let $(\freeL W,d)$ and $(\freeL W',d')$ be two dgLs satisfying~\eqref{eq:condL}.
Say that $(\freeL W', d')$ is a bigraded extension of $(\freeL W,d)$ if $W' = W \oplus \bar{W}$ as bigraded $R$-modules and $d'\!\mid_W = d$.

\begin{theorem} \label{thm:sep}
  Let $L$ be a dgL satisfying~\eqref{eq:condL}. Then $L$ has a bigraded extension $L'$ such that $L'$ is separated and the inclusion $L \injects L'$ induces an isomorphism on homology.
\end{theorem}

\begin{corollary}
  Let $X$ be a simply connected space with finite rational homology in each dimension and finite (rational) LS category. Then $X$ has a separated Lie model.
\end{corollary}

Define
\begin{equation} \label{eq:eLi}
  \eL_i = (HL_{i-1} \amalg \freeL V_{i,*}, \tilde{d}_i),
\end{equation}
where $L \amalg L'$ denotes the free product of Lie algebras (i.e., their coproduct), $\tilde{d}_i|_{HL_{i-1}}=0$ and $\tilde{d}_i|_{V_{i,*}}$ is given in \eqref{eq:tilded}.

\begin{theorem} \label{thm:HLvs}
Let $L = (\freeL V,d)$ be a separated dgL. Let 
\begin{equation*}
\hat{L}_i = \freeL((H\eL_i)_{1,*})/[\tilde{d}_{i+1}V_{i+1,*}].
\end{equation*}
Then $HL \isom \bigoplus_i \hat{L}_i$ as $\Q$-modules.
In particular, if $L$ is a Lie model for a space $X$, then its homotopy Lie algebra has the form $\pi_*(\lX) \tensor \Q \isom \oplus_i \hat{L}_i$ as rational vector spaces.
\end{theorem}

To state a more precise version of this theorem, which elucidates some of the Lie algebra structure, we need to define some more notation.

Given an increasing filtration $\ldots \subset F_{i-1}M \subset F_iM
\subset F_{i+1}M \subset \ldots$ of an $R$-module $M$, there is an
associated graded $R$-module $\gr(M) = \bigoplus_i \gr_i M$ where
$\gr_i M = F_iM / F_{i-1}M$.
If $M$ has the structure of an algebra or Lie algebra, then there is
an induced algebra or Lie algebra structure on $\gr(M)$.
If $M$ has a separate grading, then $\gr(M)$ is bigraded.
If $M$ has a differential $d$ and $dF_iM \subset F_{i-1}M$ then the filtration is called a differential filtration and there is an induced filtration on $H(M,d)$.

Recall the definition of $L_i$ and $\tilde{d}_{i+1}$ from \eqref{eq:Li} and \eqref{eq:tilded}.
There is an increasing differential filtration $\{F_k L_i\}$ on $L_i = (\freeL(V_{\leq i}), d) \isom (\freeL V_{<i} \amalg \freeL V_i, d)$ given by $F_{-1}L_i = 0$, $F_0L_i = \freeL V_{<i}$, and for $k\geq 0$, $F_{k+1} L_i = F_k L_i + [F_k L_i, V_i]$.  This induces a filtration on $HL_i$ and $HL_i/[\tilde{d}V_{i+1}]$, where we write $L/[V]$ to denote the quotient of $L$ by the Lie ideal generated by $V \subset L$. When we write $\gr(HL_i)$ and $\gr(HL_i/[\tilde{d}V_{i+1}])$ we will always be using this filtration.

There is a short exact sequence $0 \to I \xto{i} L \xto{p} A \to 0$ of Lie algebras that splits (that is, there is a Lie algebra map $j:A \to L$ such that $pj = \id_A$) if and only if $L$ is the semi-direct product of $I$ and $A$, written $L \isom A \rsdp I$. The semi-direct product is isomorphic as modules to $A \oplus I$ and the product of elements in $A$ and $I$ is given by the action of $A$ on $I$. We will show the following.

\begin{theorem} \label{thm:sepThen}
Let $L = (\freeL V, d)$ be a free dgL over $\Q$ which is separated.
Let $L_i$, $\tilde{d}_i$, $\eL_i$, $HL_i^-$ and $HL_i^+$ be as in \eqref{eq:Li}, \eqref{eq:tilded}, \eqref{eq:eLi} and just before Definition~\ref{def:separated}.
Then for all $i$, there are Lie algebra isomorphisms
\begin{gather*} 
\gr(HL_i) \isom (H\eL_i)_0 \rsdp \freeL ((H\eL_i)_1), \quad
\gr(HL_{i+1}^-) \isom HL_i^- \rsdp \hat{L}_i,\\
(H\eL_i)_0 \isom HL_{i-1} / [\tilde{d}_iV_{i}] = HL_{i-1}/HL_{i-1}^+ \isom HL_i^-, 
\end{gather*}
where $\hat{L}_i = \freeL((H\eL_i)_1)/[\tilde{d}_{i+1}V_{i+1}] = \freeL((H\eL_i)_1)/HL_i^+$, and $HL_i^-$ is a Lie subalgebra of $HL$. 
Furthermore, if we can choose a preimage $W_i \subset HL_i$ of $(H\eL_i)_1$ such that $\freeL W_i \subset HL_i$ is a Lie ideal, then $HL_i \isom HL_i^- \rsdp \freeL W_i$.
\end{theorem}

Let $*=X_0 \subset X_1 \subset \ldots \subset X_{n-1} \subset X_n = X$
be a spherical cone decomposition (see
Section~\ref{section:background}) corresponding to a separated Lie
model of $X$. In particular $X = X_{n-1} \cup_f (\bigvee_j D^{n_j+1})$
for some map $f$ from a wedge of spheres. The separated condition
implies that $f$ is \emph{free}~\cite{bubenik:freeAndSemiInert}. Let
$L_{X_i}$ denote the homotopy Lie algebra of $X_i$. Let
$i$ denote the inclusion $X_{n-1} \incl X$. Then there is an induced
map $i_{\#}: L_{X_{n-1}} \to L_X$. In general, this map is neither
injective nor surjective~\cite{halperinLemaireInert}. Let $R(L_X)$ denote the
\emph{radical} of $L_X$. That is, the sum of the solvable ideals of
$L_X$~\cite{fhjlt:radical}.

\begin{corollary} 
  Either,
  \begin{enumerate}
  \item \label{cora} $i_{\#}: L_{X_{n-1}} \to L_X$ is surjective,
  \item \label{corb} $\dim(H\eL_n)_1 = 1$, or
  \item \label{corc} $R(L_X) \subset \im(i_{\#})$ and $L_X$ contains
    a free Lie algebra on two generators.
  \end{enumerate}
\end{corollary}

\begin{rem}
  \begin{enumerate}
  \item An equivalent condition to (\ref{cora}) is the condition that $f$ is
    \emph{inert}~\cite{halperinLemaireInert}.
  \item The only example of (\ref{corb}) seems to be $\CP^n$ (see
    Example~\ref{example:cpn}).
  \item It is known that $\dim R(L_X)_{\text{even}} \leq
    n$~\cite{fhjlt:radical}. The Avramov-F{\'e}lix
    conjecture~\cite{avramovFreeLieSubalgebras,fhtCat2} states that
    either $L_X$ is finite dimensional, or $L_X$ contains a free Lie
    algebra on two generators.
  \end{enumerate}
\end{rem}

\begin{proof}
  Applying Theorem~\ref{thm:sepThen}, $\gr(L_X) \isom \im(i_{\#}) \rsdp
  \freeL ((H\eL_n)_1)$. The result follows.
\end{proof}

\begin{example} \label{example:cpn}
We show how the well-known homotopy Lie algebras of $\CP^n$ and $\CP^{\infty}$ can be easily calculated using separated Lie models.

$\CP^{\infty}$ has the Lie model $L = (\freeL \langle v_1, v_2, v_3, \ldots \rangle,d)$ where $v_k$ has dimension $2k-1$ and $dv_k = \frac{1}{2} \sum_{i+j=k}[v_i,v_j]$. The Lie subalgebra $L_i = (\freeL \langle v_1, \ldots, v_n \rangle, d)$ is a Lie model for $\CP^n$.
Then $L_1 = (\freeL \langle v_1 \rangle, 0) = \eL_1$ and $L_2 = (\freeL \langle v_1, v_2 \rangle, d) = \eL_2$ where $dv_2 = \frac{1}{2} [v_1,v_1]$. We remark that for any dgL, $L_2$ is always separated, since $HL_1^- = 0$ and with respect to $L_2$ (not $L$), $HL_2^+ = 0$. 
Now $(H\eL_2)_{0,*} = \freeL\langle v_1 \rangle / [v_1, v_1] \isom \hat{L}_1$ and $(H\eL_2)_{1,*} \isom \Q \{ u_2 \}$ where $u_2 = [v_1, v_2]$. By Theorem~\ref{thm:sepThen}, the Lie algebra structure of $HL_2$ is determined by the action of $(H\eL_2)_0$ on $\freeL ((H\eL_2)_{1,*})$ in $H\eL_2$. Since $d[v_2,v_2] = [v_1,[v_1,v_2]]$ in $L_2$, $[v_1,u_2]=0$ in $H\eL_2$. So as Lie algebras, $HL_2 \isom \freeLab \langle v_1, u_2 \rangle$, where $\freeLab$ denotes the free abelian Lie algebra.

Next we will show by induction that $L_n$ is separated, as Lie algebras, $HL_n \isom \freeLab \langle v_1, u_n \rangle$, where $u_n$ has dimension $2n$, and $HL_n^- = \freeLab \langle v_1 \rangle$. By definition $\eL_n \isom (\freeLab \langle v_1, u_{n-1} \rangle \ \amalg \ \freeL \langle v_n \rangle, \tilde{d}_n)$ where $\tilde{d}_n v_1 = \tilde{d}_n u_{n-1} = 0$ and $\tilde{d}_n v_n = u_{n-1}$. So $HL_{n-1}^- \cap HL_{n-1}^+ = \freeLab \langle v_1 \rangle \cap \freeL \langle u_{n-1} \rangle = 0$. Since $(H\eL_n)_0 \isom \freeLab \langle v_1 \rangle$, $(H\eL)_1 \isom \Q \{ [v_1, v_n] \}$ and $[v_1,[v_1,v_n]=0$ by the Jacobi identity, the claim follows.

Finally, this implies that $HL \isom \freeLab \langle v_1 \rangle$.
\end{example}

\begin{example} \label{example:productOfSpheres}
  In the next example we show that the Lie model obtained from the minimal spherical cone decomposition of a product of spheres is separated, and use this to calculate the homotopy Lie algebra for the wedges of spheres of various ``thickness''.

Let $X = \prod_{i=1}^r S^{n_i}$, where $n_i \geq 2$. Let $N = \dim X = \sum_{i=1}^r n_i$. 
Let $X_k$ denote the subcomplex of $X$ consisting of 
those points in $X$ such that at least $r-k$ of the coordinates are the basepoint.  In particular, $X_1$ is the wedge, and $X_{r-1}$ is the fat wedge. Also, $X_{k+1}$ can be obtained from $X_k$ by attaching a wedge of spheres.
Then,
\begin{equation} \label{eq:prodSpheres}
* = X_0 \subset X_1 \subset  \ldots \subset X_r =X
\end{equation}
is a spherical cone decomposition for $X$. (It is minimal, since the cone length of $X$ is bounded below by the rational LS category of $X$, which is $r$.) In addition, $X_k$ has a
Lie model $L_k = \left( \freeL \left< \oplus_{i=1}^k V_i \right>,
  d\right)$, where $V_i = \Q \{ \alpha_{i,j}\}$ with the
$\alpha_{i,j}$ in one-to-one correspondence with the $i$-fold products
of the spheres $S^{n_{\ell}}$ and the dimension of $\alpha_{i,j}$,
$|\alpha_{i,j}|$, is one less than the dimension of the corresponding
product, and $dV_i \subset \freeL \left< \oplus_{j=1}^{i-1} V_j \right>$.
Let $L_{X_k}$ denote the homotopy Lie algebra of $X_k$. Then $HL_k
\isom L_{X_k}$ and $HUL_k \isom UHL_k \isom UL_{X_k}$ where $U$
denotes the universal enveloping algebra functor.

For a non-negatively graded $\Q$-vector space $M$, let $M(z)$ be the
formal power series $\sum_{i\geq 0} (\dim M_i) z^i$. $M(z)$ is called
the Hilbert series for $M$, and $H_*(\Omega X;\Q)(z) = (UL_X)(z)$ is
called the Poincar\'{e} series for $\Omega X$.  Let $M(z)^{-1}$ denote the
power series $\frac{1}{M(z)}$. Define $A_i(z)$ by
\begin{equation*}
\prod_{i=1}^r \left(1-z^{n_i-1}x\right) = \sum_{i=0}^r A_i(z)x^i.
\end{equation*}
Let $A(z) = \prod_{i=1}^r ( 1-z^{n_i-1}) = \sum_{i=0}^r A_i(z)$, and finally define $B_k(z) = (-z)^{k-1} \sum_{i=k+1}^r A_i(z)$.

\begin{theorem} \label{thm:prodSpheres}
  As Lie algebras, $L_{X_1} \isom \freeL \langle x_1, \ldots, x_r \rangle$, where $|x_i| = n_i-1$, $L_{X} \isom \oplus_{i=1}^r \freeL \langle x_i \rangle$ and for $r\geq 3$, $L_{X_{r-1}} \isom L_X \amalg \freeL \langle u \rangle$ where $|u| = N-2$. Furthermore, for $k\geq 2$, $UL_{X_k}(z)^{-1} = A(z) - B_k(z)$.
\end{theorem}

\begin{proof}
  The first two Lie algebra isomorphisms in the statement of the
  theorem follow directly from the well-known formulas for the
  homotopy Lie algebra of a wedge and a product.

  The remainder of the statement of the theorem is obtained
  inductively.  We assume that $k\geq 2$.  Notice that the inclusion
  $X_1 \incl X$ induces a surjection $L_{X_1} \onto L_X$. So, $L_X
  \isom HL_2^- \incl HL_k^-$. Assume that $L_k$ is separated, which is
  trivial for $k=2$. Then using Theorem~\ref{thm:sepThen} one can show
  that in fact, $HL_k^- \isom L_X$. Thus $HL_k^- \incl L_X$ and hence
  $HL_k^- \cap HL_k^+ = 0$.  Therefore $L_{k+1}$ is separated. Thus,
  the minimal cone decomposition~\eqref{eq:prodSpheres} yields a
  separated Lie model for $X$.

  For $k\geq 2$, the Poincar\'{e} series for $\Omega X_k$ is obtained as
  follows. By Theorem~\ref{thmSepThenStrFree} and \cite[Lemma 3.8 and
  Theorem 3.5]{bubenik:freeAndSemiInert}) we can apply Anick's
  formula~\cite[Theorem 3.7]{anickthesis}, 

  \begin{equation*}
\begin{split}
  (UL_{X_k})(z)^{-1} &= (U(H\eL_k)_0)(z)^{-1} - [ V_{k+1}(z) +  z[(UL_{X_{k-1}})(z)^{-1} - \\ & \quad (U(H\eL_k)_0)(z)^{-1}]]\\
  &= A(z) + (-z)^{k-1}A_k(z) - z[(UL_{X_{k-1}})(z)^{-1} - A(z)],
\end{split}
\end{equation*}
where the second equality is by Theorem~\ref{thm:sepThen} and since $(UL_X)(z)^{-1} = A(z)$. By
induction, this is equal to $A(z) - B_k(z)$.

The Lie algebra isomorphism for $L_{X_{r-1}}$ follows since the cell
attachment from $X_{r-2}$ to $X_{r-1}$ is
semi-inert~\cite{bubenik:freeAndSemiInert}.
\end{proof}
We remark that the fact that the top-cell attachment of $X$ is
inert~\cite{halperinLemaireInert} is witnessed by $\eL_r =
(L_X \amalg \freeL \langle u,v \rangle, d)$, where $d|_{L_X} = 0$ and
$dv=u$.

  From Theorem~\ref{thm:prodSpheres} it is easy to check that a fat wedge of odd-dimensional spheres has the maximum possible gap in the rational homotopy groups~\cite[Theorem 33.3]{fht:rht} if and only if all the spheres have the same dimension.
\end{example}

\begin{example} \label{example:connectedSum}
  For the final example, we calculate the homotopy Lie algebra of a
  connected sum of products of spheres. For example, a Lie group $M$
  is rationally equivalent to a product of odd spheres (and so $L_M$ is a
  free abelian Lie algebra).

  For $s\geq 2$ and $1\leq i\leq s$, let $M_i$ be simply-connected, of
  dimension $N$ and rationally equivalent to a product of at least
  three spheres (e.g. $SU(n)$, $n\geq 4$).  Let $X = \#_{i=1}^s M_i$.
  Applying the previous example gives:
\begin{theorem}
  As Lie algebras,
\[
L_X \isom \coprod_{i=1}^s L_{M_i} \amalg \freeL \langle u_1, \ldots u_s \rangle / (u_1 + \ldots + u_s),
\]
where $|u_i| = N-2$.
\end{theorem}

\begin{proof}
  By the previous example $M_i$ has a separated Lie model of the form
  $(L_{(i)} \amalg \freeL\langle v_i \rangle, d)$ such that $\eL_{r_i}
  = (L_{M_i}, 0) \amalg (\freeL \langle u_i, v_i \rangle, \tilde{d})$
  with $\tilde{d} v_i = u_i$. Thus $X$ has a separated Lie model
  $(\coprod_{i=1}^s L_{(i)} \amalg \freeL \langle v \rangle, d)$ and
  $\eL_r = (\coprod_{i=1}^s L_{M_i}, 0) \amalg (\freeL \langle u_1,
  \ldots u_s, v \rangle , \tilde{d}v = u_1 + \ldots + u_s)$, where
  $r=\max_i r_i$.
\end{proof}

In particular, $L_X$ contains a free Lie algebra on $2s-1$ generators.
So $X$ satisfies the Avramov--F{\'e}lix
conjecture.
\end{example}

In the appendix we state and prove a generalization of the Schreier property of free Lie algebras: that any Lie subalgebra is also free, which may be of independent interest. This generalization is used in the proof of Theorem~\ref{thm:sepThen}.

\subsection{Acknowledgments}

{I would like to thank Kathryn Hess for many helpful discussions and Greg Lupton and John Oprea for help in rewriting the introduction and suggesting Examples \ref{example:productOfSpheres} and \ref{example:connectedSum}.}

\section{Background} \label{section:background}

In their landmark papers \cite{quillen:rht, sullivan:infinitesimalComputationsInTopology}, Quillen and Sullivan construct algebraic models for rational homotopy theory.  Sullivan constructs a contravariant functor $A_{PL}$ which serves as a fundamental bridge between topology and algebra. For a space $X$, $A_{PL}(X)$ is a commutative cochain algebra which has the property that $H(A_{PL}(X)) \isom H^*(X;\Q)$ as algebras. Quillen gives a construction for a differential graded Lie algebra (dgL) $(\freeL V, d)$, where $\freeL V$ denotes the free Lie algebra on a rational vector space $V$, one of whose properties is that $H(\freeL V,d) \isom \pi_*(\lX) \tensor \Q$. We will follow convention and call this a free dgL even though it is almost always not a free object in the category of \emph{differential} graded Lie algebras. 
For an excellent reference on these models and their applications, the reader is referred to~\cite[Parts II and IV]{fht:rht}.

A \emph{quasi-isomorphism} is a morphism which induces an isomorphism in homology. 
A \emph{Lie model} for a space $X$ is a differential graded Lie algebra $(L,d)$  equipped with a quasi-isomorphism $m: C^*(L,d) \qi A_{PL}(X)$. Here, $C^*(L,d) = \Hom(C_*(L,d),\Q)$, is the contravariant functor induced by Quillen's functor $C_*$. $C_*(L,d)$ is called the Cartan-Eilenberg-Chevalley construction on $(L,d)$. Quillen's free dgL above is a Lie model. Given a Lie model $(L,d)$ for a space $X$ and a quasi-isomorphism $(L,d) \qi (L',d')$, their is an induced quasi-isomorphism $C^*(L',d') \qi C^*(L,d) \qi A_{PL}(X)$. In particular, a quasi-isomorphic bigraded extension of a Lie model for $X$ is also a Lie model for $X$.

It is the Samelson product on $\pi_*(\lX)$, which corresponds the the Whitehead product under the canonical isomorphism $\pi_*(\lX) \isom \pi_{*+1}(X)$, which gives it the structure of a graded Lie algebra, which we call the \emph{(rational) homotopy Lie algebra}.

A \emph{rational homotopy equivalence} is a continuous map $f:X \to Y$ such that $\pi_*(f) \tensor \Q$ is an isomorphism. Two spaces $X$ and $Y$ are said to have the same \emph{rational homotopy type}, written $X \rht Y$, if they are connected by a sequence of rational homotopy equivalences (in either direction). The \emph{Lusternik Schnirelmann (LS) category} of a space $X$, denoted $\cat(X)$, is the smallest integer $n$ such that $X$ is the union of $n+1$ open subsets, each contractible in $X$. The \emph{rational LS category} of a space $X$, denoted $\catQ(X)$, is the smallest integer $n$ such that there exists a space $Y$ with $\cat(Y) = n$ and $X \rht Y$. So $\catQ(X) \leq \cat(X)$.
A space $X_n$ is said to be a \emph{spherical $n$-cone} if there exists a
sequence of spaces  
\begin{equation} \label{eqnXN}
* = X_0 \subset X_1 \subset \ldots \subset X_n
\end{equation}
such that for $k=0\ldots n-1$, $X_{k+1} = X_k \cup_{f_{k+1}} \left( \bigvee_j D_{j,k}^{n_{j}+1} \right)$, where $D^n$ denotes the $n$-dimensional disk and $f_{k+1}:\bigvee_j S^{n_j}_{j,k} \to X_k$ is an attaching map.
A space $X$ is said to have \emph{rational cone length} $n$, written
$\clQ(X) = n$, if $n$ is the smallest number such that $X \rht X_n$
for some spherical $n$-cone $X_n$.
(This is equivalent to the more usual definition of rational cone length, see~\cite[Proposition 28.3]{fht:rht}, for example.)
It is a theorem of Cornea~\cite{cornea:coneLengthAndLScat} that if $\catQ(X) = n$ then $\clQ(X) = n$ or $n+1$.
A CW complex is said to have finite type if it has finitely many cells in each dimension.

A free dgL $(\freeL V,d)$ is said to have \emph{length} $N$ if we can decompose $V_* = \oplus_{i=1}^N V_{i,*}$ such that $dV_{i,*} \subset \freeL (\oplus_{j=1}^{i-1} V_{j,*})$, where $dV_{1,*}=0$. Note that $V$ is now bigraded, though $(\freeL V,d)$ is typically not a bigraded dgL. We will call the first gradation \emph{degree} and the second, usual gradation \emph{dimension}.
For any free dgL $(\freeL V,d)$ of length $N$ there is a spherical $N$-cone $X$ such that $(\freeL V, d)$ is a Quillen model for $X$. We call $(\freeL V,d)$ the cellular Lie model of $X$. Furthermore, any spherical $N$-cone has a cellular Lie model of length $N$. We outline the construction as follows. 
For each $n \geq 1$ let the $(n+1)$-cells $D^{n+1}_{\alpha}$ of $X$ correspond to a basis $\{v_{\alpha}\}$ of $V_n$. Let $X_n$ denote the $n$-skeleton of $X$. By induction, $(\freeL V_{<n}, d)$ is a cellular Lie model of $X_n$. The attaching maps $f_{\alpha}: S^n \to X_n$ are given by representatives of $[dv_{\alpha}] \in H_{n-1}(\freeL V_{<n}, d) \isom \pi_{n-1}(\lX_n) \isom \pi_{n}(X_n)$.

A dgL is said to be of \emph{finite type} if it is finite in each dimension. It is said to be \emph{connected} if it is concentrated in strictly positive dimensions. Given two (Lie) algebras $L$ and $L'$, let $L \amalg L'$ denote their free product (i.e., coproduct). We remark that $\freeL V \amalg \freeL V' \isom \freeL (V \oplus V')$.

\section{Replacing dgLs with separated dgLs}

To prove Theorem~\ref{thm:sep}, we will need the following lemmas. Let $\kk$ be a field of characteristic $0$.

\begin{lemma} 
Let $W$ be a free $\kk$-module and let $\alpha,\beta \in (\freeL W, d)$
with $d\beta = \alpha$. 
Let $W' = W \oplus \kk\{a,b\}$ where $|a|=|\beta|$ and $|b|=|\beta|+1$.
Define $(\freeL W', d')$ by letting $d'$ be an extension of $d$
determined by taking $d'a = \alpha$ and $d'b= a - \beta$. 
Then the inclusion $(\freeL W, d) \injects (\freeL W', d')$ is a quasi-isomorphism.
\end{lemma}

\begin{proof}
Consider the following commutative diagram.
\begin{equation*}
  \xymatrix{ (\freeL W, d) \ar[rr]^{\phi} \ar[dr]^{\psi} & & (\freeL (W \oplus     \kk\{a,b\}), d') \\
    & (\freeL (W \oplus \kk\{\hat{a},\hat{b}\}),\hat{d}) \ar[ur]^{\theta} &
  }
\end{equation*}
where $\hat{d}\hat{b} = \hat{a}$, $\phi$ and $\psi$ are injections and $\theta$ is defined on generators as follows: $\theta|_W = \id_W$, $\theta \hat{a} = a-\beta$, and $\theta \hat{b} = b$.
It is easy to check that $\theta$ is a chain map and a dgL isomorphism. Furthermore, 
\[
(\freeL (W \oplus \kk\{\hat{a},\hat{b}\}),\hat{d}) \isom 
(\freeL W, d) \amalg (\freeL<\hat{a},\hat{b}>,\hat{d}) \simeq (\freeL W, d).
\]
Thus $\psi$ is a quasi-isomorphism. It follows that $\phi$ is one as well.
\end{proof}

\begin{corollary} \label{corExtension}
Let $L = (\freeL W, d)$ with $\{\alpha_j, \beta_j\}_{j \in J} \subset
\freeL W$ where $d\alpha_j = 0$ and $d\beta_j=\alpha_j$.
Taking $\bar{W} = \kk\{a_j, b_j\}_{j \in J}$ with $|a_j| = |\beta_j|$ and
$|b_j| = |\beta_j|+1$, 
let $L' = \freeL(W \oplus \bar{W}, d')$ where $d'|_W = d$, $d'a_j =
\alpha_j$, and $d'b_j = a_j - \beta_j$.
Then $L' \simeq L$.
\end{corollary}

Given an $\kk$-module $M$, let $ZM$, $BM$ denote the $\kk$-submodules of
cycles and boundaries.

\begin{lemma} \label{lemmaNAttach}
Let $L = (\freeL W,d)$ and let $V = \{v_j\}_{j \in J} \subset HL_{n}$
with $v_j \neq 0$.
For each $v_j$ choose a representative cycle $\hat{v}_j \in ZL$.
Let $L' = (\freeL W',d')$ where $W' = W \oplus \kk\{a_j,b_j\}_{j \in
J}$, $d'\!\mid_W = d$, $da_j = \hat{v}_j$ and $b_j$ is in dimension $n+2$.
Then $H_{\leq n}L' \isom H_{\leq n}L/V$.
\end{lemma}

\begin{proof}
Since $Z_{\leq n}L' = Z_{\leq n}L$ and $B_{\leq n}L' \isom
B_{\leq n}L \oplus \kk\{da_j\}_{j \in J}$, $H_{\leq n}L' \isom
H_{\leq n}L / V$.
\end{proof}

\begin{lemma} \label{lemmaVS}
Given $\kk$-modules $C \subset A \subset B$ and $D \subset B$, let $p: B
\to B/C$ denote the quotient map.
If $D \cap A \subseteq C$ then $pD \cap pA = 0$.
\end{lemma}

\begin{proof}
Let $x \in pD \cap pA$. 
Then there are $y \in D$, and $y' \in A$ such that $py = py' = x$.
Let $z = y - y'$.
Since $pz = 0$, $z \in C \subset A$.
Thus $y = y' + z \in A \cap D \subseteq C$.
Therefore $x = py = 0$.
\end{proof}

The proof of theorem~\ref{thm:sep} will rely on a two-step inductive
procedure given in Proposition~\ref{propSeparated}.
To help with the book-keeping, we introduce the following definitions.
We will use the notation $HL_i^+$ and $HL_i^-$ defined just before Definition~\ref{def:separated}.

\begin{defn}
  Let $(\freeL W, d)$ be a dgL satisfying~\eqref{eq:condL}.  Say that $(\freeL W, d)$ is \emph{$k$-separated} if for all $i<k, \ HL_i^+ \cap HL_i^- =   0$.  Say $(\freeL W, d)$ is \emph{$(k,n)$-separated} if it is $k$-separated and   $H_{<n}L_k^+ \cap H_{<n}L_k^- = 0$.
\end{defn}

Let $\delta_{i,j} = 1$ if $i=j$ and let it be $0$ otherwise.
Let $s$ denote the suspension homomorphism. 
That is, $(sM)_k = M_{k-1}$.

\begin{lemma} \label{lemmaSeparated}
Let $L = (\freeL W, d)$ be a dgL of length $N$ that is $(k,n)$-separated. 
Let $V$ be the component of $HL_k^+ \cap HL_k^-$ in dimension $n$.

Then there exists a bigraded extension $\hat{L} = (\freeL \hat{W},
\hat{d}) \supset L$ of length $N+\delta_{k+2,N+1}$ with 
$\hat{W}_i = W_i$ for $i \neq k, k+2$, 
$\hat{W}_k = W_k \oplus \bar{W}$,
$\hat{W}_{k+2} = W_{k+2} \oplus s\bar{W}$,
where $\bar{W}$ is in dimension $n+1$,
such that $L' \simeq L$.
Furthermore
$H_{\leq n}\hat{L}_k \isom H_{ \leq n}L_k / V$, and
\begin{equation} \label{eqnLv0} 
H_{\leq n}\hat{L}_k^+ \cap H_{\leq n}\hat{L}_k^- = 0.
\end{equation}
\end{lemma}

\begin{proof}
Let $\{\alpha_j\}_{j \in J}$ be a basis (as an $\kk$-module) for $V$.
Since $\alpha_j \in HL_k^-$, $\alpha_j$ has a preimage
$\alpha'_j \in HL_{k-1}$.
Let $\alpha''_j$ be a representative cycle in $L_{k-1}$ for
$\alpha'_j$.
Since $\alpha_j \in HL_k^+$, there is a $\beta_j \in
L_{k+1}$ such that $d\beta_j = \alpha''_j$.
Let $\{a_j, b_j \}_{j \in J}$ be pairs of elements of bidegree $(k,n+1)$ and
$(k+2,n+2)$ respectively.
Define $\hat{L} = (\freeL \hat{W}, \hat{d}) \supset L$ by letting 
$\hat{W}_i = W_i$ for $i \neq k, k+2$, $\hat{W}_k = W_k \oplus
\kk\{a_j\}$ and $\hat{W}_{k+2} = W_{k+2} \oplus \kk\{b_j\}$.
Extend $d$ to $\hat{d}$ by letting $da_j = \alpha''_j$ and $db_j = a_j -
\beta_j$. 
Note that $\hat{L}$ has length $N+\delta_{k+2,N+1}$.
By Corollary~\ref{corExtension}, $\hat{L} \simeq L$.

By Lemma~\ref{lemmaNAttach}, 
\begin{equation} \label{eqnLxn}
(H_{\leq n}\hat{L}_k) \isom (H_{\leq n}L_k) / V.
\end{equation}
In other words, in dimensions $\leq n$ the map from $HL_k$ to
$H\hat{L}_k$ is just the quotient by $V$.
Thus, $H_{\leq n}\hat{L}_k^+ \isom H_{\leq n}L_k^+ / V$, and $H_{\leq n}\hat{L}_k^- \isom H_{\leq n}L_k^- / V$.
Applying Lemma~\ref{lemmaVS} to 
$V \subset H_{\leq n}L_k^- \subset H_{\leq n}L_k$
and $H_{\leq n}L_k^+ \subset H_{\leq n}L_k$, one gets that $H_{\leq n}L_k^+ / V \cap H_{\leq n}L_k^- / V = 0$.
\end{proof}

\begin{defn}
Say that a bigraded $\kk$-module $M$ is in the \emph{$(k,n)$-region} if
$W_{*,\leq n} = W_{\geq k+3, *} = W_{k+1, n+1} = W_{k+2,n+1} = 0$.
Let $L = (\freeL W, d)$ be a free dgL.
Say that $L' = (\freeL W',d')$ is a \emph{$(k,n)$-extension} of $L$ if
$W' = W \oplus \bar{W}$ as bigraded modules, where $\bar{W}$ is in the $(k,n)$-region, and $d'|_W = d$.
\end{defn}

\begin{prop} \label{propSeparated}
Let $L = (\freeL W, d)$ of length $N$ that is $(k,n)$-separated.
Let $V$ be the component of 
$HL_k^+ \cap HL_k^-$ in dimension $n$.

(a) Then there is an $(k,n)$-extension of $L$, 
$L' = (\freeL W',d')$ of length $N+\delta_{k+2,N+1}$ which is 
$(k,n+1)$-separated, such that $L' \simeq L$.
Furthermore $H_{\leq n}L'_k \isom H_{\leq n}L_k / V$.

(b) In addition there an $(k,n)$-extension of $L$, 
$L'' = (\freeL W'', d'')$ of length $N+\delta_{k+2,N+1}$ which is 
$(k+1)$-separated, such that $L'' \simeq L$.
Furthermore $H_{\leq n}L''_k \isom H_{\leq n}L_k / V$, and if 
$H_{\leq n}L_{k+1}^+ \cap H_{\leq n}L_{k+1}^- = 0$ then
$H_{\leq n}{L''}_{k+1}^+ \cap H_{\leq n}{L''}_{k+1}^- = 0$.
\end{prop}

\begin{rem}
Recall that $d(W_{k+1}) \subset \freeL(W_{\leq k})$ and that
$HL_k^+$ is the ideal generated by the
images of $W_{k+1}$ by the induced map $W_{k+1} \to HL_k$.
All of the elements in this ideal have preimages in $\freeL (W_{\leq k})$.
$V$ consists of those dimension $n$ elements in this ideal which have
preimages in lower filtration.
\end{rem} 

\begin{proof}
We prove the proposition by induction on $k$. 
The statement of the proposition is trivial if $k=0$.
We will assume the statement of the proposition is true for $k-1$.

Let $L$, $V$ be as in the statement of the proposition.
By Lemma~\ref{lemmaSeparated}, there exists $\hat{L} = (\freeL \hat{W},
\hat{d}) \supset L$ of length $N+\delta_{k+2,N+1}$ with 
$\hat{W}_i = W_i$ for $i \neq k, k+2$, 
$\hat{W}_k = W_k \oplus \bar{W}$,
$\hat{W}_{k+2} = W_{k+2} \oplus s\bar{W}$,
where $\bar{W}$ is in dimension $n+1$,
such that $L' \simeq L$.
Furthermore
$H_{\leq n}\hat{L}_k \isom H_{\leq n}L_k / V$, and
\begin{equation} \label{eqnZeroInHat}
H_{\leq n}\hat{L}_k^+ \cap H_{\leq n}\hat{L}_k^- = 0.
\end{equation}

Since $\hat{L}_{k-1} = L_{k-1}$ and $L$ is $k$-separated, $\hat{L}$ is
(at least) $(k-1)$-separated. 
Also $H\hat{L}_{k-2} = HL_{k-2}$ and $H\hat{L}_{k-1} = HL_{k-1}$.

Let $\hat{\mathcal L} = H\hat{L}_{k-1}^+ \cap H\hat{L}_{k-1}^-$. Since $L$ is $k$-separated, this equals $[\hat{d} \bar{W}] \cap HL_{k-1}^-$.
Since $\bar{W}$ is $n$-connected, $\mathcal{\hat{L}}$ is
$(n-1)$-connected, but it is not necessarily $n$-connected.
In the non-trivial case, there are elements of $V$ which not only have
preimages in filtration $k-1$ but also have preimages in filtration
$k-2$.
Let $\hat{V} = \mathcal{\hat{L}}_n$. 

By induction there exists a $(k-1,n)$-extension $L' = (\freeL V',
d') \supset \hat{L}$ of length $N + \delta_{k+2,N+1}$ that is
$k$-separated such that $L' \simeq \hat{L}$,
and since $H_{\leq n}\hat{L}_k^+ \cap H_{\leq n}\hat{L}_k^- = 0$, 
$H_{\leq n}{L'}_k^+ \cap H_{\leq n}{L'}_k^- = 0$.
So $L'$ is a $(k,n)$-extension of $L$ such that $L' \simeq L$.
Furthermore $H_{\leq n}L'_{k+1} \isom H_{\leq n}\hat{L}_{k+1} /
\hat{V} = H_{\leq n}L_{k-1} / \hat{V}$. 
This proves part (a) of the statement.

To prove part (b) of the statement we simply iterate part (a).
By iterating (a) we get a sequence of dgLs
\[ L = L^{(0)} \subset L^{(1)} \subset \ldots \subset L^{(i)} \subset
L^{(i+1)} \subset \ldots
\]
where $L^{(i)}$ is $(k,n+i)$-separated, 
$L^{(i+1)}$ is a $(k,n+i)$-extension of $L^{(i)}$,
and $L^{(i+1)} \simeq L^{(i)}$.
Furthermore, $H_{\leq n}L^{(i)}_k \isom H_{\leq n}L_k / V$, and if 
$H_{\leq n}L_{k+1}^+ \cap H_{\leq n}L_{k+1}^- = 0$ then $H_{\leq n}{L^{(i)}}_{k+1}^+ \cap H_{\leq n}{L^{(i)}}_{k+1}^- = 0$.
Let $L'' = \cup_i L^{(i)}$.
Then $L''$ is a $(k+1)$-separated $(k,n)$-extension of $L$ and $L'' \simeq L$.
Furthermore, $H_{\leq n}L^{(i)}_k \isom H_{\leq n}L_k / V$, and if 
$H_{\leq n}L_{k+1}^+ \cap H_{\leq n}L_{k+1}^- = 0$ then $H_{\leq n}{L''}_{k+1}^+ \cap H_{\leq n}{L''}_{k+1}^- = 0$.
\end{proof}

\begin{proof}[Proof of Theorem~\ref{thm:sep}]
Since any dgL is $1$-separated,
Theorem~\ref{thm:sep} follows by applying
Proposition~\ref{propSeparated}, $N$ times.
\end{proof}

\section{Properties of separated dgLs}

In this section we will use Theorems \ref{thm:sep} and
\ref{thmSchreier} to prove Theorem~\ref{thm:sepThen}.
We will defer to the appendix the proof of Theorem \ref{thmSchreier}, which is a generalization the well-known result that a Lie subalgebra of a free Lie algebra is also a free Lie algebra. We will use \eqref{eq:Li}, \eqref{eq:tilded}, \eqref{eq:eLi}, and the notation $HL_i^+$ and $HL_i^-$ defined just before Definition~\ref{def:separated}.

\begin{defn}
Let $(\freeL W,d)$ be a free dgL over a field.
Say that $(\freeL W,d)$ is \emph{strongly free} if for all $i$,
$HL_i^+ \subset L_{M_i}$ is a free Lie subalgebra. 
\end{defn}

\begin{theorem} \label{thmSepThenStrFree}
Let $L$ be a free dgL over $\Q$ which is separated.
Then $L$ is strongly free and for all $i$, there are Lie algebra isomorphisms
\begin{equation*} 
\gr(HL_i) \isom (H\eL_i)_0 \rsdp \freeL ((H\eL_i)_1) \text{, and }
(H\eL_i)_0 \isom HL_{i-1} / [\tilde{d}W_{i}] \isom HL_i^-.
\end{equation*}
\end{theorem}

We will prove Theorem~\ref{thmSepThenStrFree} by induction.
A main part of the induction will use the following theorem, which is
a special case of the main algebraic result from~\cite{bubenik:freeAndSemiInert}.
Note that $U$ denotes the universal enveloping algebra functor.

\begin{theorem}[{\cite[Theorem 3.12]{bubenik:freeAndSemiInert}}] \label{thmHUL}
Over the field $\Q$, let $L' = (L \amalg \freeL V_1, d)$, where $dL \subset L$, $dV_1 \subset L$ and $HUL \isom UL_0$, such that there is an induced map
$d':V_1 \to L_0$.
Let $\eL = (L_0 \amalg \freeL V_1, d')$.
If $[d'V_1] \subset L_0$ is a free Lie algebra,
then as algebras
\[
\gr(HUL') \isom UH\eL \text{, and }
H\eL \isom (H\eL)_0 \rsdp \freeL (H\eL)_1 \text{ and } (H\eL)_0 \isom L_0/[d' V_1]
\]
as Lie algebras.  
\end{theorem}
The associated graded structure on $HUL'$ is that induced from the filtration on $L'$ given by $F_{-1}L'=0$, $F_0L' = L$ and $F_{i+1}L' = F_iL' + [V_1,F_iL']$.

\begin{proof}[Proof of Theorem~\ref{thmSepThenStrFree}]
Let $L = (\freeL W, d)$ be a free dgL over $\Q$ which is separated.
Recall that $\tilde{d}_i: W_i \to HL_{i-1}$ and $\eL_i =
(HL_{i-1} \amalg \freeL W_i, \tilde{d}_i)$.

Assume that $L_{i}$ is strongly free, $\gr(HL_i) \isom (H\eL_i)_0 \rsdp
\freeL ((H\eL_i)_1)$, and $(H\eL_i)_0 \isom HL_{i-1} / [\tilde{d}W_i]
\isom HL_i^-$.
This is trivial for $i \leq 1$.

Note that $L_{i+1} = (L_i \amalg \freeL W_{i+1}, d)$.
Both $L_{i+1}$ and $UL_{i+1}$
can be filtered by the `length in $W_{i+1}$' filtration.
That is, let $F_{-1}(L_{i+1}) = 0$, $F_0 L_{i+1} = L_i$, and for
$i\geq 0$, $F_{i+1}(L_{i+1}) = F_i(L_{i+1}) + [F_i(L_{i+1}), W_{i+1}]$. 
This induces a similar filtration on $UL_{i+1}$.

By assumption, $\gr(HL_i) \isom HL_i^- \rsdp \freeL ((H\eL_i)_1)$.
Since $L_{i+1}$ is separated, $HL_i^+ \cap
HL_i^- = 0$.
Thus by Theorem~\ref{thmSchreier}, $HL_i^+$ is a free Lie
algebra.
Hence $L_{i+1}$ is strongly free.

Then by Theorem~\ref{thmHUL}, as algebras
\begin{gather*}
\gr(HUL_{i+1}) \isom U \left( (H\eL_{i+1})_0 \rsdp \freeL ((H\eL_{i+1})_1) \right) \text{, and }\\ 
(H\eL_{i+1})_0 \isom HL_i / [\tilde{d}_{i+1} W_{i+1}].
\end{gather*}
 
It is a result of Quillen's that $HUL_{i+1} \isom UHL_{i+1}$.
Thus $\gr(HUL_{i+1}) \isom \gr(UHL_{i+1}) \isom U\gr(HL_{i+1})$.
From this it follows that $U\gr(HL_{i+1}) \isom U\left( (H\eL_{i+1})_0 \rsdp \freeL
((H\eL_{i+1})_1) \right)$.

For any Lie algebra $L$, $UL$ has a canonical cocommutative Hopf
algebra structure.
Let $P$ denote the primitive functor.
Over $\Q$, the composition $PU$ is the identity
functor~\cite{milnorMooreHopfAlgebras}. 
Therefore $\gr(HL_{i+1}) \isom (H\eL_{i+1})_0 \rsdp \freeL
((H\eL_{i+1})_1)$.

It follows that $HL_{i+1}^- \isom (H\eL_{i+1})_0$.
This finishes the inductive step.
\end{proof}


We will now use Theorem \ref{thmSepThenStrFree} to prove Theorem~\ref{thm:sepThen}.

\begin{proof}[Proof of Theorem~\ref{thm:sepThen}]
Let $L=(\freeL V,d)$ that is separated.
Recall that $L_i = \freeL (V_{\leq i, *})$, 
$\tilde{d}:V_i \xto{d} ZL_{i-1} \onto HL_{i-1}$, and $\eL_i =
(HL_i \amalg \freeL V_i, \tilde{d})$.

Let $N_i = (H\eL_i)_1$ and recall that $HL_i$ is filtered by the length in $V_i$ filtration.  By Theorem~\ref{thmSepThenStrFree}, $(\freeL V,d)$ is strongly free and for all $i$, it satisfies $\gr(HL_i) \isom HL_i^- \rsdp \freeL N_i$ where $HL_i^- \isom HL_{i-1} / {[\tilde{d}_iV_i]}$.

Also by Theorem~\ref{thmSepThenStrFree}, $HL_{i+1}^- \isom
HL_i / [\tilde{d}_{i+1}V_{i+1}]$, 
which has a filtration induced by the filtration on $HL_i$.
Thus we have the following short exact sequence of filtered,
graded Lie algebras: 

\[
0 \to [\tilde{d}_{i+1}V_{i+1}] \to HL_i \to HL_{i+1}^- \to 0,
\]
which induces a short exact sequence of bigraded Lie algebras:
\[
0 \to \gr( [\tilde{d}_{i+1}V_{i+1}] ) \to \gr(HL_i) \to \gr(
HL_{i+1}^- ) \to 0.
\]
Recall that 
$\gr(HL_i) \isom HL_i^- \rsdp \freeL N_i$.

Since $(\freeL V,d)$ is separated, it follows that $[\tilde{d}_{i+1}V_{i+1}] \cap HL_i^- =
0$, and thus $\gr([\tilde{d}_{i+1}V_{i+1}]) \subset \freeL N_i$ as Lie algebras.
Since any Lie subalgebra of a free Lie algebra is automatically free,
$\gr([\tilde{d}_{i+1}V_{i+1}]) \isom \freeL K_{i+1}$ for some $\Q$-module
$K_{i+1} \subset \freeL N_i$.
Therefore $\gr(HL_{i+1}^-) \isom HL_i^- \rsdp (\freeL
N_i / \freeL K_{i+1})$ as Lie algebras.
Let $\hat{L}_i = \freeL N_i / \freeL K_{i+1}$.

By Theorem~\ref{thmSepThenStrFree} there is a split short exact sequence of Lie algebras
\[
0 \to \freeL ((H\eL_i)_1) \to \gr(HL_i) \to HL_i^- \to 0
\]
Let $W_i$ be a preimage of $(H\eL_i)_1$ in $HL_i$.
We have a short exact sequence of modules
\[
0 \to \freeL W_i \to HL_i \to HL_i^- \to 0
\]
but this is not necessarily a short exact sequence of Lie algebra since $\freeL W_i$ may not be a Lie ideal. Equivalently, the projection $HL_i \to HL_i^-$ may not be a Lie algebra morphism. 
However, if $\freeL W_i$ is a Lie ideal in $HL_i$ then $HL_i \isom HL_i^- \rsdp \freeL W_i$.
\end{proof}

\appendix

\section{A generalized Schreier property}

In this chapter we give a simple criterion which we prove guarantees
that certain Lie subalgebras are free. 

It is a well-known fact that any (graded) Lie subalgebra of a (graded)
Lie algebra is a free Lie algebra
\cite{mikhalevSubalgebras, shirshovSubalgebras, mikhalevZolotykh}
This is often referred to as the \emph{Schreier property}.
In this chapter we generalize this result to the following.

\begin{theorem} \label{thmSchreier}
Over a field $\F$, let $L$ be a finite-type graded Lie algebra with filtration
$\{F_k L\}$ such that $\gr(L) \isom L_0 \rsdp \freeL V_1$ as Lie
algebras, where $L_0 = F_0 L$ and $V_1  = F_1 L / F_0 L$.
Let $J \subset L$ be a Lie subalgebra such that $J \cap F_0 L = 0$.
Then $J$ is a free Lie algebra.
\end{theorem}

Before proving this theorem,  we prove the following lemma.

\begin{lemma} \label{lemmaSchreier}
Let $J$ be a finite-type filtered Lie algebra such that $\gr(J)$ is a
free Lie algebra.
Then $J$ is a free Lie algebra.
\end{lemma}

\begin{proof}
By assumption there is an $\F$-module $\bar{W}$ such that $\gr(J)
\isom \freeL \bar{W}$.

Let $\{\bar{w}_i\}_{i \in I} \subset \gr(J)$ be an $\F$-module basis
for $\bar{W}$.
Let $m_i = \deg(\bar{w}_i)$.
That is, $\bar{w}_i \in F_{m_i}J / F_{m_i-1} L$.
For each $\bar{w}_i$ choose a representative $w_i \in F_{m_i} J$.
Let $W = \F\{w_i\}_{i \in I} \subset J$.

Then there is a canonical map $\phi: \freeL W \to J$.
Grade $\freeL W$ by letting $w_i \in W$ be in degree $m_i$.
Then $\phi$ is a map of filtered objects and there is an induced map
$\theta: \freeL W \to \gr(J)$.
However the composite map
\[ \freeL W \xto{\theta} \gr(J) \isomto \freeL \bar{W}
\]
is just the canonical isomorphism $\freeL W \isomto \freeL \bar{W}$.
So $\theta$ is an isomorphism.
Therefore $\phi$ is an isomorphism and $J$ is a free Lie algebra.
\end{proof}

\begin{proof}[Proof of Theorem~\ref{thmSchreier}]
The filtration on $L$ filters $J$ by letting 
\[ F_k J = J \cap F_k L.
\]
From this definition it follows that the inclusion $J \incl L$ induces
an inclusion $\gr(J) \incl \gr(L)$.
So $\gr(J) \incl \gr(L) \isom L_0 \rsdp \freeL V_1$.
Since $J \cap F_0 L = 0$ it follows that $(\gr J)_0 = 0$ and $\gr(J)
\incl (\gr J)_{\geq 1} \isom \freeL V_1$.  
By the Schreier property $\gr(J)$ is a free Lie algebra.
Thus by Lemma~\ref{lemmaSchreier}, $J$ is a free Lie algebra.
\end{proof}

The following corollary is a special case of this theorem.

\begin{corollary} \label{corSchreier}
Over a field $\F$, if $J \subset L_0 \rsdp \freeL (L_1)$ is a Lie
subalgebra such that $J \cap L_0 = 0$ then $J$ is a free Lie algebra.\
\end{corollary}

Note that since $J$ is not necessarily homogeneous with respect to
degree $J \cap L_0 = 0$ does not imply that $J \subset \freeL L_1$. 


\newcommand{\etalchar}[1]{$^{#1}$}

\end{document}